\documentclass[11pt]{article}

\usepackage{latexsym}
\usepackage{amsfonts}
\usepackage{amssymb}

\newtheorem{theorem}{Theorem}[section]
\newtheorem{proposition}{Proposition}[section]

\newcommand{\bP}{{\bf P}}
\newcommand{\bE}{{\bf E}}

\newcommand{\tildebar}[1]{\stackrel{\simeq}{#1}}
\newcommand{\dbar}[1]{\stackrel{=}{#1}}
\newcommand{\Pg}{{\bf P}^{\bar{g}^*}}
\newcommand{\Eg}{{\bf E}^{\bar{g}^*}}
\newcommand{\Varg}{\mbox{{\bf Var}}^{\bar{g}^*}}
\newcommand{\sfrac}[2]{{\textstyle \frac{#1}{#2}}}
\newcommand{\err}{\max \left(
    \sfrac{1}{n},\sqrt{\sfrac{\hat{\tau}}{nm}} \right)} 

\title{How to Combine Fast Heuristic Markov Chain Monte Carlo
with Slow Exact Sampling}
\author{Antar Bandyopadhyay \\ and \\ David J. Aldous\thanks{This material is based upon work supported by the National Science Foundation under Grant No. 9970901} \\
\\
\\
       University of California\\
       Department of Statistics\\
        367 Evans Hall \# 3860\\
       Berkeley CA 94720-3860}

\begin{document}
\maketitle

\vspace{0.1in}

\begin{abstract}

Given a probability law $\pi$ on a set $S$ and a function
$g : S \rightarrow R$, suppose one wants to estimate the mean
$\bar{g} = \int g \ d\pi$. The Markov Chain Monte Carlo method
consists of inventing and simulating a Markov chain with stationary
distribution $\pi$. Typically one has no \emph{a priori} bounds on the chain's
mixing time, so even if simulations suggest rapid mixing one cannot
infer rigorous confidence intervals for $\bar{g}$.
But suppose there is also a separate method which (slowly) gives samples
exactly from $\pi$.
From $n$ exact samples, one could
immediately get a confidence interval of length
$O(n^{-1/2})$. But one can do better.  Use each exact sample as the initial
state of a Markov chain, and run each of these $n$ chains for $m$
steps. We show how to construct confidence intervals which are always
valid, and which, if the (unknown) relaxation time of the chain 
is sufficiently small relative to $m/n$, have length
$O(n^{-1} \log n)$ with high probability.

\end{abstract}


\newpage
\section{Background}
\label{INT}

Let $\pi$ be a given probability distribution on a set
$S$. Given a function $g : S \rightarrow R$, we want to estimate
its mean $\bar{g} := \int_S g(s) \pi(ds)$. 
As we learn in elementary statistics,
one can obtain an estimate for $\bar{g}$ by taking
samples from $\pi$ and using the sample average $g$-value as an estimator. 
But algorithms
which sample exactly from $\pi$ may be prohibitively slow.
This is the setting for the
{\em Markov chain Monte Carlo} (MCMC) method,
classical in statistical physics and over the last ten years studied
extensively as statistical methodology \cite{CSI00,MCMCIP,liu-MCMC,RC99}.
In MCMC one designs a Markov chain on state-space $S$ to have
stationary
distribution $\pi$.  Then the sample average $g$-value over a long run of the
chain is a heuristic estimator of $\bar{g}$.
Diagnostics for assessing length of run required,
and expressions for heuristic confidence intervals, form a
substantial part of MCMC methodology \cite{DiscMCMC}.
In general one cannot make such estimates rigorous, because one cannot
eliminate the possibility that all the samples seen come from some small
part of the state space which is almost disconnected from the remainder.
Rigorous estimates typically require an a priori bound on
some notion of the chain's \emph{mixing time} (e.g. the relaxation time 
defined at (\ref{def-relax})); and while there
is now substantial theoretical literature on mixing times
\cite{meRWG,behrends,DSC97a} it deals with settings more tractable than
most statistical applications.

This paper investigates the interface between rigor and heuristics
in a particular (perhaps artificial) context.
Suppose we have a guess $\hat{\tau}$ for the mixing time of the
chain, based on simulation diagnostics or heuristic estimates
\cite{Garren00} or some non-rigorous
mathematical argument.
Suppose we have some separate scheme (an {\em exact sampler}) 
which gives independent samples
exactly from $\pi$.
Imagining that sampling $\hat{\tau}$
steps of the chain is roughly equivalent to sampling once from $\pi$,
it is natural to
consider the ratio
\[ \rho =
\frac{
\mbox{cost of one exact sample}}{
\mbox{cost of $\hat{\tau}$ steps of the chain}} \]
where {\em cost} refers to computational time.
If $\rho < 1$ then one would just use the exact sampler and forget
MCMC.
If $\rho$ is extremely large then we might not be able to afford even
one exact sample, and we are forced to rely on MCMC
(this is the setting typically envisaged in MCMC).
This paper addresses the remaining context, where $\rho$ is large
but not extremely large; in other words, we can afford
to simulate many steps of the chain (enough to make estimates
heuristically good) but can afford only a few exact samples.
In the case of sampling from general $d$-dimensional densities,
for instance,
exact samplers (e.g. based on rejection sampling using some
tractable comparison density) are typically feasible only for
small $d$, 
and MCMC is used for large $d$, so there should always be a
window of $d$-values which fits our ``$\rho$ large, but not extremely large"
context.

In this context, we could just
use the exact sampler to get $n$ independent samples from $\pi$.
Then the sample average $g$-value provides an estimate
of $\bar{g}$ with $O(n^{-1/2})$ error.
But instead, suppose we use these $n$ independent samples as
initial states and generate $n$ independent $m$-step realizations
of the Markov chain. 
If diagnostic tests suggest that mixing occurs in
$\hat{\tau}$ steps then we have an ``effective sample size"
of $(n \times m/\hat{\tau})$ and the heuristic estimate
of error (when we use the overall sample average $g$-value as 
an estimator) would be 
$O(\sqrt{\hat{\tau}/(nm)})$.
Our main result, Theorem \ref{T1}, shows that in a certain sense
such error bounds can be made rigorous.

\section{Results} 
The discussion in section \ref{INT} provides conceptual context for our result,
but let us now state the (rather minimal) mathematical assumptions
for the result.
For simplicity we assume the state space $S$ is finite
(though since our results are non-asymptotic they must extend to the
general case without essential change).
We assume (for reasons explained in Section \ref{dis-T1}) the 
function $g:S \to R$ is bounded, so by rescaling we may assume
\begin{equation}
0 \leq g(\cdot) \leq 1 .
\label{g-bounded}
\end{equation}
We assume the Markov chain is reversible, that is
to say its transition matrix $K$ satisfies
\begin{equation}
 \pi_i k_{ij} = \pi_j k_{ji}, \quad \forall i,j  .\label{def-reversible}
\end{equation}
These are the only background assumptions for 
{\em validity} of the conservative confidence interval
given at (\ref{valid}).
That is, there are no ``implicit asymptotics"; and we do not even need
to assume $K$ is irreducible.
The {\em length} of the confidence interval will depend on
the data, i.e. the realizations of the chain, but (\ref{I-length}) bounds the length in terms of
the {\em relaxation time}
of the chain, defined as
\begin{equation}
\tau_2 := (1-\lambda_2)^{-1} \label{def-relax} \end{equation}  where 
$1= \lambda_1 \geq \lambda_2 \geq \cdots \geq \lambda_s \geq -1$ 
are the eigenvalues of $K$,
and $\tau_2 < \infty$ if $K$ is irreducible.
\begin{theorem}
\label{T1}
Assume (\ref{g-bounded},\ref{def-reversible}).
Take $n \geq 3$, $m \geq 1$ and $0<\alpha<1$ and $\hat{\tau} \geq 1$.
Based on $2n$ exact samples from $\pi$ and $2mn$ steps
of the $K$-chain, we can construct an interval $I$ such that
\begin{equation}
\bP \left( \bar{g} \not\in I \right) \leq \alpha
\label{valid}
\end{equation}
and
\[ \bP \left( \mbox{length}(I) > k_{\alpha} \, \err \, \log n \right) \]
\begin{equation}
\leq
3 n (n+1) \exp \left( - \frac{m}{48 n \tau_2} \log^2 n \right),
\label{I-length}
\end{equation}
where $k_{\alpha} := 2 \left( \sqrt{2/\alpha} + \log (4/\alpha) \right)$.
\end{theorem}

\subsection{Discussion}
\label{dis-T1}
To interpret Theorem \ref{T1}, suppose we take 
$m = n \hat{\tau}$, so that we use $2n^2 \hat{\tau}$ steps of the chain. 
Then the ``target length"(in the left side of (\ref{I-length})) of our 
confidence interval will be $O(n^{-1} \log n)$.
The theorem guarantees a confidence interval that is always valid,
and guarantees that, if $\tau_2$ is indeed not more than $\hat{\tau}$,
then the length of the confidence interval
will likely not exceed the target length. 
This contrasts with the $O(n^{-1/2})$ length confidence interval
obtained by using the exact sampler only, and gets close to the $O(n^{-1})$ length of the
heuristic confidence interval.

Admittedly our numerical error bounds are too crude to be much use in
practice. For example in order to make the bound in (\ref{I-length})
smaller than say $5\%$, one would have to generate at least $300$
exact samples, which is most often not practicable. So the result is
primarily of theoretical interest, in particular because of its
similarity to the idea \cite{RS94} of self-testing algorithms. That
paper describes an algorithm for generating random self-avoiding
walks. As the authors write \cite{RS94} ``While there are a number of
Monte Carlo algorithms used to solve these problems in practice, these
are heuristic and their correctness relies on unproven conjectures. In
contrast, our algorithm is shown rigorously to produce answers with
specific accuracy and confidence. Only the efficiency of the algorithm
relies on a widely believed conjecture, and a novel feature is that
this conjecture can be \emph{tested} as the algorithm proceeds.'' In our MCMC
setting, we cannot estimate rigorously the actual value of $\tau_2$, but we can
self-justify inferences based on estimated $\tau_2$.

On a more technical note, let us outline why Theorem \ref{T1}
gives close to the best possible bounds on confidence interval length.
Indeed, we claim that the best one could hope for is length
of order
\begin{equation}
\max \left( \frac{1}{n}, \sqrt{\frac{\tau_2}{nm}} \right).
\label{LBlength}
\end{equation}
The point is that there are two different ``obstacles" to sharp
estimation.
First, consider the eigenvector $g_2$
associated with eigenvalue $\lambda_2$.
It is easy to estimate the variance of the overall sample average
when $g = g_2$, and this variance works out as
order $\frac{\tau_2}{nm}$;
so we should not hope to have smaller estimation error than the
corresponding standard deviation $\sqrt{\frac{\tau_2}{mn}}$.
Second, suppose some subset $A$ of the state space, 
with $\pi(A) = 1/n$, is almost disconnected.
Then it is not unlikely that all $n$ exact samples, and hence
the $n$ realizations of the chains, miss $A$, and so
a contribution $E_\pi [g(\cdot) 1_A]$ to $\bar{g}$ would be
``invisible" to our simulations,
and so this contribution is an unavoidable source of possible
error when using sample averages as estimators.
Our assumption (\ref{g-bounded}) that $g$ is bounded was intended
as the simplest way of bounding this error -- bounding it as
order $1/n$.

So Theorem \ref{T1} shows that, if our initial guess $\hat{\tau}$ is indeed
roughly close to $\tau_2$, then our rigorous confidence interval's length will be 
roughly of the minimal order (\ref{LBlength}), 
up to $\log n$ terms.
In Section \ref{application} we give a natural ``adaptive" variation
in which we prescribe two numbers 
$\hat{\tau} < \hat{\tau}_{{\rm max}}$,
where as before $\hat{\tau}$ is a heuristic estimate of $\tau_2$,
and where $2n^2 \hat{\tau}_{{\rm max}}$ is the maximum number of steps of the chain
that we would be willing to simulate.
Theorem \ref{T2} gives an always-valid confidence interval which, if
$\tau_2$ is indeed small relative to $\hat{\tau}_{{\rm max}}$,
will have length of order $n^{-1} \log n$ and will require order 
$n^2 \tau_2$ steps of the chain.

\subsection{Outline of construction and proof}
Recall the ``procedure" of simulating $n$ realizations of $m$ steps
of the Markov chain, starting from $n$ exact samples from $\pi$.
The construction of the confidence interval $I$ in Theorem \ref{T1}
can be summarized as follows.\\
(i) Perform this procedure once, and find the overall average
$g$-value -- call it  $\bar{g}^*$.\\
(ii) Perform the procedure again, and for $1 \leq i \leq n$
let $A_i$ be the average $g$-value over the $i$'th $m$-step
realization. \\ (iii)  Test whether 
$|A_i - \bar{g}^*| \leq \frac{\log n}{\sqrt{r(n,m)}}$
for every $i$, where $r(n,m) := \min \left( n, \frac{m}{\hat{\tau}} \right)$.
If so, report a ``short" confidence interval
$\left[ {\rm ave}_i A_i \pm O\left(\err \, \log n \right) \right]$;
if not, report a ``long" confidence interval
$[ {\rm ave}_i A_i \pm O(\frac{1}{\sqrt{n}}) ]$.

\vspace{0.08in}
\noindent
To analyze the validity of the confidence interval,
the key point is that after observing the event
``$|A_i - \bar{g}^*| \leq \frac{\log n}{\sqrt{r(n,m)}}$"
happening $n$ times out of $n$, we can be confident that its
probability is $1 - O(1/n)$.
This allows us to truncate $A_i$ at $\bar{g}^* \pm 
\frac{\log n}{\sqrt{r(n,m)}}$, and then the sample average of $n$ truncated
variables has s.d. of order
$\frac{\log n}{\sqrt{r(n,m)}} \times \frac{1}{\sqrt{n}} = \err \,
\log n $.

Finally, to bound the chance of not reporting the short confidence
interval we need to bound the chance of a truncation 
being needed, and a bound can be derived from large deviation estimates
for reversible chains.

\subsection{An adaptive version}
\label{application}
In the procedure underlying Theorem \ref{T1} we make a single guess
$\hat{\tau}$ and hope that the ``good event" which leads to a short
confidence interval will happen; if not, we settle for a long confidence
interval.
A natural variation is to specify that, if the ``good event" fails,
then repeat the process with $\hat{\tau}$ replaced
by $2\hat{\tau},4 \hat{\tau}, 8 \hat{\tau},\ldots$ and continue until
the ``good event" happens or until we reach some predetermined limit
on numbers of steps of the chain.

\begin{theorem}
\label{T2}
Assume (\ref{g-bounded},\ref{def-reversible}). Take $n \geq 5$, 
$ 0 < \alpha < 1$, and $1 \leq \hat{\tau} \leq \hat{\tau}_{{\rm max}} 
= 2^a \hat{\tau}$, where $a \geq 0$ is an integer. Then based on $2n$
exact samples from $\pi$, and $2n^2 \times M$ steps of the $K$-chain,
where $M$ is a random variable taking values in
$\{ \hat{\tau}, 2 \hat{\tau}, 2^2 \hat{\tau}, \ldots, 2^a \hat{\tau} \}$, 
we can define an interval $I$, such that 
\begin{equation} 
\bP \left( \bar{g} \not \in I \right) \leq \alpha, 
\label{con-1}
\end{equation} 
and 
\begin{equation} 
\mbox{length}(I) \leq k_{\alpha}^a \, \frac{\log n}{n}
\mbox{\ \ \ \ if \ \ \ } 1 \leq \frac{M}{\hat{\tau}} < 2^a, 
\label{con-2}
\end{equation} 
and 
\begin{equation} 
\bP \left( M > 96 \, \tau_2 \vee \hat{\tau} \right) 
\leq 3 n (n+1) \exp \left( - \log^2 n \right).
\label{con-3}
\end{equation} 
where $k_{\alpha}^a = 2 \left( \sqrt{2(a+1)/\alpha} 
+ \log (4(a+1)/\alpha) \right)$. 
\end{theorem}

So we are prescribing the maximum number of steps of the chain to be
$2n^2 \hat{\tau}_{{\rm max}}$.
The bound in (\ref{con-3}) is less than $0.05$ for $n=8$ and
goes to zero very rapidly as $n$ increases.  So if
$\hat{\tau}_{{\rm max}}$ is indeed large compared to $\tau_2$ then by
generating a small number of exact samples one can construct a
confidence interval for $\bar{g}$ which will be ``short" with high probability, 
and the number of steps of the Markov chain required
will be $O(n^2 ( \tau_2 \vee \hat{\tau}))$.

\section{Proof of Theorem \ref{T1}}
\label{Proof-T1}

\subsection{Construction of the confidence interval}
\label{sec-constr}

Let $\{Z_1^*, Z_2^*, \ldots, Z_n^*\}$ be $n$ independent samples from
$\pi$. For $1 \leq i \leq n$ let $(X_{ij}^*)_{j=0}^{m-1}$ be a
reversible Markov chain with initial state $X_{i0}^* = Z_i^*$; these $n$ Markov
chains are independent. Define 
$$A_i^* := \frac{1}{m} \sum_{j=0}^{m-1} g(X_{ij}^*), \mbox{\ \ } 1 \le
i \le n,$$ 
and
\begin{equation}
\bar{g}^* := \frac{1}{n} \sum_{i=1}^n A_i^* .
\label{defg*}
\end{equation} 
$\bar{g}^*$ is our initial guess for $\bar{g}$. 

Now re-run the entire simulation independently to get  
$\{Z_1, Z_2, \ldots, Z_n\}$, another set of $n$ independent samples
from $\pi$, and $(X_{ij})_{j=0}^{m-1}$ another independent but
identically distributed family of $n$ reversible Markov chains each
with initial state $X_{i0} = Z_i$. We further define 
$$A_i := \frac{1}{m} \sum_{j=0}^{m-1} g(X_{ij}), \mbox{\ \ } 1 \le i \le n,$$
and 
\begin{equation} 
\dbar{A} := \frac{1}{n} \sum_{i=1}^n A_i.
\label{Adbar}
\end{equation}
Truncate each $A_i$ to get 
\begin{equation} 
\tilde{A}_i := \left\{ \begin{array}{cl}
                          A_i & \mbox{\ \ \ if \ \ \ } \vert A_i - \bar{g}^*
                          \vert \le \frac{\log n}{\sqrt{r(n,m)}} \\
                          \bar{g}^* & \mbox{\ \ \ otherwise}
                          \end{array} \right. 
\label{Atildes}
\end{equation}
where $r(n,m) := \min \left( n, \frac{m}{\hat{\tau}} \right)$.
Let 
\begin{equation} 
\tildebar{A} := \frac{1}{n} \sum_{i=1}^{n} \tilde{A}_i . 
\label{Atildebar}
\end{equation} 
Write
$N_n = \sum_{i=1}^{n} I ( A_i \neq \tilde{A}_i ) $
for the number of truncations, and call
the event $G_n := [ N_n = 0 ]$ the \emph{good event}.  

Define $h : (\mathbb{N} \cup \{0\} ) \times \mathbb{N} \times (0,1)
\to [0, \infty)$ by
\begin{equation}
h(z,n;\alpha) := \left\{ \begin{array}{ccc}
                            {\displaystyle \frac{z}{n} +
                            \frac{c_{\alpha}}{\sqrt{n}}} & \mbox{if} &
                            z \neq 0 \\ 
                            {\displaystyle \frac{d_{\alpha}}{n}} &
                              \mbox{if} & z = 0 
                            \end{array} \right.,
\label{defh}
\end{equation} 
where 
$$c_{\alpha} := \frac{1}{\sqrt{2 \alpha}} \mbox{\ \ and \ \ } 
d_{\alpha} := \log \frac{2}{\alpha}. $$
In the next section we shall prove
\begin{proposition}
\label{P1}
For any $b>0$
\begin{equation} 
\bP \left[ \vert \tildebar{A} - \bar{g} \vert > b \,
  \err \, \log n + h(N_n,n;\alpha) \right] \leq \frac{1}{b^2} + \alpha. 
\label{T2.1}
\end{equation} 
\end{proposition}
Replacing $\alpha$ by $\alpha/2$ in Proposition \ref{P1}
and setting 
$b=\sqrt{2/\alpha}$, we see that the confidence interval 
\begin{equation} 
I:= \tildebar{A} \pm \left( \sqrt{\frac{2}{\alpha}} \err \, \log n + 
   h(N_n,n;\alpha /2) \right) 
\label{CI}
\end{equation} 
satisfies the requirement of ({\ref{valid}) that
$\bP(\bar{g} \notin I) \leq \alpha$.
If $N_n=0$ then
the length of this confidence interval is
\begin{equation}
 2 \left( \sqrt{\frac{2}{\alpha}} \, \err \, \log n + 
\frac{d_{\alpha / 2}}{n} \right) .
\label{length}
\end{equation}
Notice that (\ref{length}) is bounded by 
$k_{\alpha} \, \err \, \log n$;   
here we use assumption $n \geq 3$ which implies $\log n > 1$.
So to prove (\ref{I-length}) and complete the proof of Theorem \ref{T1},
it is enough to prove (in section \ref{sec-LD})
\begin{proposition}
\label{P2}
\[ \bP(N_n >0) \leq
3 n (n+1) \exp \left( - \frac{m}{48 n \tau_2} \log^2 n \right). \]
\end{proposition}

\subsection{Proof of Proposition \ref{P1}}

We denote the conditional probability, conditional expectation and
conditional variance given $\bar{g}^*$ by $\Pg$, $\Eg$ and $\Varg$
respectively. 

Observe that under $\Pg$, the random variables
$\tilde{A}_1, \tilde{A}_2, \ldots, \tilde{A}_n$ are i.i.d. Thus 
$\Eg(\tildebar{A}) = \Eg (\tilde{A}_1)$,
and $\Varg (\tildebar{A}) = \frac{1}{n} \Varg (\tilde{A}_1)$. 
But $\Varg (\tilde{A}_1) = \Varg ( \tilde{A}_1 - \bar{g}^*) \leq
      \left( \frac{\log n}{\sqrt{r(n,m)}} \right)^2$, because
$ \vert \tilde{A}_1 - \bar{g}^* \vert \leq \frac{ \log n }{ \sqrt{r(n,m)}}$. 
So by Chebyshev's inequality, we get
$$ \Pg \left[ \vert \tildebar{A} - \Eg ( \tilde{A}_1 )
              \vert  > b \, \err \, \log n \right] \leq \frac{1}{b^2}, $$
and by taking expectation 
\begin{equation} 
\bP \left[ \vert \tildebar{A} - \Eg ( \tilde{A}_1 )\vert >
  b \, \err \, \log n \right] \leq \frac{1}{b^2}.
\label{Cbound}
\end{equation} 

Now we want to estimate 
$ \vert \Eg ( \tilde{A}_1 ) - \bar{g} \vert $. 
From the definitions (\ref{Adbar}) and (\ref{Atildebar}),
$$ ( \dbar{A} - \tildebar{A} ) = \frac{1}{n} \sum_{i=1}^n A_i I(A_i
\ne \tilde{A}_i ) - \frac{1}{n} \sum_{i=1}^n \bar{g}^* I(A_i \ne
\tilde{A}_i).$$ 
Since $g$ takes values in $[0,1]$, 
\begin{equation} 
- \frac{1}{n} N_n \mbox{\ } \le \mbox{\ }  \Bigl( \dbar{A} -
\tildebar{A} \Bigr) \mbox{\ } \le \mbox{\ } \frac{1}{n} N_n.
\label{temp}
\end{equation} 
Now $\dbar{A}$ is independent of $\bar{g}^*$, thus taking
conditional expectation given $\bar{g}^*$ in (\ref{temp}) we get 
\begin{equation} 
\left\vert \Eg(\tilde{A}_1) - \bar{g} \right\vert \mbox{\ } \leq
\mbox{\ } p_n(\bar{g}^*),
\label{Ebound}
\end{equation} 
where
$p_n(\bar{g}^*) := \Pg \Bigl( \vert A_1 - \bar{g}^* \vert > \log n /
\sqrt{r(n,m)} \Bigr) $. 

Under $\Pg$ we have $ \mbox{\ } N_n \sim 
\mbox{{\bf Binomial}}\left( n , p_n({\bar{g}^*})\right)$, and hence
\begin{eqnarray} 
\Pg\left( p_n(\bar{g}^*) > \frac{d_{\alpha}}{n}, N_n = 0 \right)
&=& \Bigl( 1 - p_n(\bar{g}^* ) \Bigr)^n 1_{(p_n(\bar{g}^*) > d_\alpha /n)} 
\nonumber\\ 
&\leq &\left( 1 - \frac{d_{\alpha}}{n} \right)^n 1_{(d_{\alpha}/n \leq
  1)} \nonumber\\
&\leq & e^{-d_\alpha} \mbox{\ \ since \ } ( 1 - x ) \leq e^{-x}
\mbox{\ \ } \forall \mbox{\ } 0 \leq x \leq 1 \nonumber\\
&=& \frac{\alpha}{2} \mbox{\ \ by definition of \ } d_{\alpha}.
\label{t1}
\end{eqnarray}
Further, \begin{eqnarray}
\lefteqn{           \Pg \left( p_n(\bar{g}^*) >
           \frac{N_n}{n}+\frac{c_{\alpha}}{\sqrt{n}}, N_n > 0 \right)
}\\
             & \le & \Pg \Biggl( \frac{N_n}{n} +
           \frac{c_{\alpha}}{\sqrt{n}} < p_n(\bar{g}^*) \Biggr)
           \nonumber \\
             & \leq & \Pg\left( \left\vert N_n - n p_n(\bar{g}^*)
           \right\vert > \sqrt{n} c_{\alpha} \right) \nonumber \\
             & \leq & \frac{p_n(\bar{g}^*) ( 1 - p_n(\bar{g}^*))}
           {c_{\alpha}^2} \mbox{ by Chebyshev's inequality} \nonumber \\
             & \leq & \frac{1}{4 c_{\alpha}^2} \nonumber \\
             & = & \frac{\alpha}{2} \mbox{\ \ by definition of \ }
           c_{\alpha}. \label{t2} 
           \end{eqnarray} 
Taking expectations of the conditional probabilities in (\ref{t1}) and
(\ref{t2}) we get 

$$ \bP \left( p_n(\bar{g}^*) > \frac{d_{\alpha}}{n}, N_n = 0 \right) 
\leq \frac{\alpha}{2} $$
and
$$ \bP \left( p_n(\bar{g}^*) > \frac{N_n}{n} +
  \frac{c_{\alpha}}{\sqrt{n}}, N_n > 0 \right) \leq \frac{\alpha}{2}. $$
Thus by definition of $h(\cdot)$
\begin{equation} 
\bP \Bigl( p_n(\bar{g}^*) > h(N_n,n;\alpha) \Bigr) \leq \alpha .
\label{pbound}
\end{equation} 
And hence from (\ref{Cbound}), (\ref{Ebound}) and (\ref{pbound}) we get 
\begin{eqnarray*}
&      & \bP \left( \vert \tildebar{A} - \bar{g} \vert 
              > b \,\err \, \log n + h(N_n,n;\alpha) \right) \\ 
& \leq & \bP \left( \vert \tildebar{A} - \bE( \tilde{A}_1 \vert
         \bar{g}^* ) \vert > b \, \err \, \log n \right) \\ 
&      & \qquad \qquad \mbox{} + \bP \Bigl( p_n(\bar{g}^*) > h(N_n, n;
         \alpha) \Bigr) \\  
& \leq & \frac{1}{b^2} + \alpha .
\end{eqnarray*}

\hspace{4.5in} $\square$

\vspace{.1in}

\subsection{Proof of Proposition \ref{P2}}
\label{sec-LD}
Clearly
\begin{eqnarray}
\bP(N_n>0)&\leq& n \bP \left( |A_1-\bar{g}^*|> \frac{\log
  n}{\sqrt{r(n,m)}} \right) \nonumber \\  
&\leq & n \left[ \bP \left( \vert A_1 - \bar{g} \vert > \frac{\log
  n}{2 \sqrt{r(n,m)}} \right) \right. \nonumber \\
&     & \qquad \qquad + \left. \bP \left( \vert \bar{g}^* - \bar{g} \vert >
        \frac{\log n}{2 \sqrt{r(n,m)}} \right) \right].
\label{LD1}
\end{eqnarray}


To bound the terms of (\ref{LD1}) we use a large deviation
bound for sample averages of reversible Markov chains.
Lezaud \cite{lezaud98} equation (2)
gives a one-sided bound for 
$A_1 = m^{-1} \sum_{j=0}^{m-1} g(X_{1j})$:
$$ \bP \Bigl(  A_1 - \bar{g}  > \lambda  \Bigr) 
\leq  \exp \left( \frac{1}{5 \tau_2}  - 
\frac{\lambda^2 m }{12 \tau_2} \right)  , \quad \lambda > 0 .
$$
Since $\tau_2 \geq 1/2$ always, and $2e^{2/5} < 3$,
we deduce the two-sided bound
\begin{equation}
 \bP \Bigl(  \vert A_1 - \bar{g} \vert  > \lambda  \Bigr) 
\leq 3  \exp \left( - \frac{\lambda^2 m }{12 \tau_2} \right)  , \quad
\lambda > 0  . \label{LLD}
\end{equation}
So in particular
\begin{eqnarray} 
\bP \left(  \vert A_1 - \bar{g} \vert  > \frac{\log n}{2 \sqrt{r(n,m)}}
\right) & \leq & 3 \exp \left( - \frac{m}{48 n \tau_2} \,
  \frac{n}{r(n,m)} \log^2 n \right) \nonumber \\
       & \leq & 3  \exp \left( -\frac{ m }{48 n \tau_2} \log^2 n  \right) .
\label{LD2}
\end{eqnarray}

Also, for $\lambda > 0$, 
\begin{eqnarray*}
\bP \left( \vert \bar{g}^* - \bar{g} \vert > \lambda \right) 
&  \leq & n \bP \left( \vert A_1^* - \bar{g} \vert > \lambda \right) \\
 &  =   & n \bP \left( \vert A_1 - \bar{g} \vert > \lambda \right) \\
 & \leq & 3 n \exp \left( - \frac{\lambda^2 m}{12 \tau_2} \right) \mbox{ by }
 (\ref{LLD}). 
\end{eqnarray*}
So in particular
\begin{eqnarray}
\bP \left( \vert \bar{g}^* - \bar{g} \vert > \frac{\log n}{2 \sqrt{r(n,m)}}
\right) & \leq & 3 n \exp \left( - \frac{m}{48 n \tau_2} \, \frac{n}{r(n,m)}
                \log^2 n \right) \nonumber \\
       & \leq & 3 n \exp \left( - \frac{m}{48 n \tau_2} \log^2 n
                \right). \label{LD3}
\end{eqnarray}
Substituting (\ref{LD2}) and (\ref{LD3}) into (\ref{LD1})
gives the bound asserted in Proposition \ref{P2}.

\hspace{4.5in} $\square$

\section{Proof of Theorem \ref{T2}}
For the first part of the procedure for constructing the
confidence interval $I$, simulate
$\{Z_i^*, 1 \leq i \leq n\}$ and $\{Z_i, 1 \leq i \leq n\}$
as at the start of section \ref{sec-constr}.  
This part of the procedure is not repeated.
Then for $k \in \{0,1,2,\ldots,a\}$, let 
$\{ X_{ij}^* \, \vert \, 1 \leq i \leq n, \, 1 \leq j \leq m_k := 2^k n
\hat{\tau} \}$ be realizations of chains started at $X_{i0}^*=Z_i^*$ ; and 
let $\{ X_{ij} \, \vert \, 1 \leq i \leq n, \, 1 \leq j \leq m_k := 2^k n
\hat{\tau} \}$ be realizations of chains started at $X_{i0}=Z_i$; each 
simulated until time $m_k := 2^k n \hat{\tau}$. 
Repeat definitions from Section \ref{Proof-T1}: for $k = 0,1,\ldots a$ 
define $A_i^{(k)}$, $\tilde{A}_i^{(k)}$, and $\tildebar{A}^{(k)}$ 
as $A_i$, $\tilde{A}_i$, and $\tildebar{A}$ respectively with 
$m = m_k$.  Note that for such $m$ we have $r(n,m) = n$. 

Let $N_n^{(k)} = \sum_{i=1}^n I ( A_i^{(k)} \not = \tilde{A}_i^{(k)} )$. 
Define
\begin{equation}  
I \left( n, k ; \alpha \right)  :=
\left[ \tildebar{A}^{(k)} \pm \left( \sqrt{\frac{2}{\alpha}} \, 
\frac{\log n}{n} + h ( N_n^{(k)}, n ; \alpha / 2 ) \right) \right],
\label{defIk}
\end{equation} 
where $h(\cdot)$ is as defined in (\ref{defh}).
This is the interval $I$ defined at (\ref{CI}) which featured in Theorem
\ref{T1},
and so by (\ref{valid}) we get that for $0 \leq k \leq a$, 
$0 < \alpha < 1$, and $n \geq 5$, 
\begin{equation} 
\bP \left( \bar{g} \not \in I(n,k;\alpha) \right) \leq \alpha. 
\label{valid-2}
\end{equation} 

Define $T := \min \{ \, 0 \leq k \leq a \, \vert \, N_n^{(k)} = 0 \, \}$, 
where we write $ T = a $ if the set is empty. 
Then define
\[M := 2^T \hat{\tau} \in 
\{ \hat{\tau}, 2 \hat{\tau}, 2^2 \hat{\tau}, \ldots, 2^a \hat{\tau} \} \]. 
\[ I := I \Bigl( n, T ; \alpha/(a+1) \Bigr). \] 
If $ 0 \leq T < a$ then $N_n^{(T)} = 0$, and hence
from (\ref{defIk}) and (\ref{defh}) we get that 
\begin{equation} 
\mbox{length}(I) \leq k_{\alpha}^a \, \frac{\log n}{n},
\end{equation} 
where $k_{\alpha}^a = 2 \left( \sqrt{2(a+1)/\alpha} 
+ \log (4(a+1)/\alpha) \right)$, so (\ref{con-2}) is satisfied. 
Further, 
\begin{eqnarray*}
\bP \Bigl( \bar{g} \not \in I \Bigr)
 & =    & \sum_{k=0}^a \bP \Bigl( \bar{g} \not \in I, T=k \Bigr) \\
 & =    & \sum_{k=0}^a \bP \Bigl( \bar{g} \not \in 
          I ( n, k; \alpha / (a+1) ), T=k \Bigr) \\
 & \leq & \sum_{k=0}^a \bP \Bigl( \bar{g} \not \in
          I \left( n, k; \alpha / (a+1) \right) \Bigr) \\
 & \leq & \sum_{k=0}^a \frac{\alpha}{a+1} = \alpha. 
\end{eqnarray*}
So (\ref{con-1}) is satisfied also. 

Now to complete the proof we observe that,
writing $96 \, \hat{\tau} \vee \tau_2$
for $96 (\hat{\tau} \vee \tau_2)$,
\begin{eqnarray} 
\bP \left( M > 96 \, \hat{\tau} \vee \tau_2 \right)
 & =    & \bP \left( T > \left\lfloor \log_2  
          \frac{96 \, \tau_2 \vee \hat{\tau}}{\hat{\tau}} 
          \right\rfloor \right) \nonumber \\
 & \leq & \bP \left( N_n^{ (\lfloor \log_2 \frac{96 \, \tau_2 
          \vee \hat{\tau}}{\hat{\tau}} \rfloor ) } > 0 \right), 
          \label{temp-1}
\end{eqnarray}
where $\lfloor x \rfloor$ denotes the greatest integer less than or equal to $x$. 

\noindent
Applying Proposition \ref{P2} with 
$m = n \times 2^{\lfloor \log_2 \frac{96 \, \tau_2 \vee
\hat{\tau}}{\hat{\tau}}\rfloor } \hat{\tau}$, we get 
\begin{eqnarray} 
\bP \left( N_n^{ (\lfloor \log_2 \frac{96 \, \tau_2 
          \vee \hat{\tau}}{\hat{\tau}} \rfloor ) } > 0 \right)  
 & \leq & 3 n (n+1) \exp \left( - \frac{2^{ \lfloor \log_2 \frac{96 \hat{\tau} \vee
          \tau_2}{\hat{\tau}} \rfloor } \, \hat{\tau}}{48 \tau_2} \log^2 n \right)
          \nonumber \\
 & \leq & 3 n (n+1) \exp \left( - \log^2 n \right) \label{temp-2}. 
\end{eqnarray}

\noindent 
The bound asserted in (\ref{con-3}) follows from (\ref{temp-1}) and 
(\ref{temp-2}). 
The number of chain steps used equals
$2 n^2 \times 2^T \hat{\tau} = 2 n^2 \times M $. 

\hspace{4.5in} $\square$

\end{document}